\documentclass[11pt]{article}
\usepackage{amsmath}
\usepackage{bm}
\usepackage{url}
\usepackage{algpseudocode}
\usepackage{hyperref}
\usepackage{graphicx}
\usepackage{authblk}
\usepackage{jlcode}
\usepackage{geometry}
\hypersetup{
    colorlinks=true, 
    linkcolor=blue,
    citecolor=blue,
    urlcolor=blue,
}

\begin{document}

\title{RationalFunctionApproximation.jl: Rational Approximation On Discrete and Continuous Domains}

\author[1]{Tobin A. Driscoll}
\affil[1]{Department of Mathematical Sciences, University of Delaware}

\hypersetup{
pdftitle = {RationalFunctionApproximation.jl: Rational Approximation On Discrete and Continuous Domains},
pdfauthor = {Tobin A. Driscoll},
pdfkeywords = {Julia, Approximation, Complex variables, Rational functions},
}

\maketitle

\begin{abstract}
Unlike polynomials, rational functions can represent functions having poles or branch cuts with root-exponential convergence and no Runge phenomenon. Recent developments of the AAA and greedy Thiele algorithms have sparked renewed interest in computational rational approximation. The \textsf{RationalFunctionApproximation} package supplies the fastest known implementations of these methods and the only arbitrary-precision ones. Combined with the \textsf{ComplexRegions} package, it can produce compact and accurate representations of a huge variety of functions over intervals, circles, or other domains in the complex plane.
\end{abstract}

\section{Introduction}

Approximation of functions is a foundational technology for scientific computing. As of this writing, the Julia ecosystem has excellent options for approximation in packages such as \textsf{Interpolations} \cite{kittisopikulInterpolationsjl2025} for piecewise polynomials, \textsf{ApproxFun}~\cite{olverPracticalFramework2014} for global polynomials, \textsf{KernelInterpolation}~\cite{lampertKernelInterpolationjlMultivariate2024} for radial basis functions, and many others.

A long-studied class of approximating functions are the rational functions, which are ratios of polynomials. Unlike polynomials, they can represent functions having branch cuts with root-exponential convergence~\cite{newmanRationalApproximation1964, herremansResolutionSingularities2023} and have no Runge requirement for crowding interpolation nodes near boundaries~\cite{platteImpossibilityFast2011}. Rational functions have been used in applications such as nonlinear eigenvalue problems, reduced-order modeling, and the solution of PDEs~\cite{antoulasInterpolatoryMethods2020,budisaRationalApproximation2022,gopalSolvingLaplace2019,hautSolvingBurgers2013,lietaertAutomaticRational2022,rodriguezPAAAAlgorithm2022,trefethenNumericalConformal2020}.

An important connection between rational and polynomial interpolants is their common representation in the barycentric formula~\cite{berrutBarycentricLagrange2004,schneiderNewAspects1986}, which allows for fast, stable evaluation. A major recent development in computational rational approximation is the AAA algorithm~\cite{nakatsukasaAAAAlgorithm2018}. This approach chooses the barycentric weights to minimize least-squares error over a test set and uses an iterative greedy selection of new interpolation nodes. The method has inspired many extensions and variations~\cite{filipRationalMinimax2018,hochmanFastAAAFast2017,nakatsukasaAlgorithmReal2020,wilberDatadrivenAlgorithms2022}. Of particular interest to this work is the extension of AAA to adaptive node selection in continuous domains in order to avert under-resolution near singularities~\cite{driscollAAARational2024}.

Rational functions have an alternative classical representation in terms of continued fractions, a fact used by Thiele at least 100 years ago to derive an interpolation scheme via inverse differences~\cite{milne-thompsonCalculusFinite1933}. Because of the connection to divided differences, Thiele's method has long been known to be vulnerable to numerical instability, depending strongly on the selection and ordering of interpolation nodes. But Salazar~\cite{celisAdaptiveThiele2023,celisNumericalContinued2024} has recently shown that a greedy iterative approach to selecting nodes can yield stable Thiele approximations. Adding a single node to a Thiele interpolant on $n$ nodes requires just $O(n)$ work, compared to $O(n^3)$ for AAA; the possibility of a faster method is intriguing, although the stability situation is not yet fully understood.

The \textsf{RationalFunctionApproximation} package implements three rational approximation methods: the AAA algorithm, the greedy Thiele algorithm, and linear least-squares approximation using prescribed poles. It builds on the \textsf{ComplexRegions} package to provide a convenient interface for approximating functions on complex domains, such as the unit disk or the exterior of a polygon. The package is compatible with extended-precision arithmetic as provided, e.g., by \textsf{BigFloat} or \textsf{MultiFloats}. The package is intended for use both as a library for applications or for research into rational approximation methods.

This paper is written about and using version 0.2.4. 

\section{Mathematical background}
\label{sec:math}

\subsection{AAA}
\label{sec:aaa}

The barycentric representation of a rational interpolant is
\begin{equation}
    \label{eq:barycentric}
    r(z) = \frac{ \; \displaystyle \sum_{j=1}^n \frac{y_j\, w_j}{z - z_j} \; }{ \displaystyle \sum_{j=1}^n \frac{w_j}{z - z_j} },
\end{equation}
where $z_1, \ldots, z_n$ are distinct interpolation nodes, $y_1, \ldots, y_n$ are the values of $r$ at those nodes, and $w_1, \ldots, w_n$ are the barycentric weights. Generically, the function $r$ is of degree $n$ in both the numerator and denominator. For a particular choice of the weights, however, $r$ is actually a polynomial of degree at most $n-1$~\cite{berrutBarycentricLagrange2004}. We can instead use the weights to improve the global quality of the approximation of $r$ to a function $f$ with values $y_j = f(z_j)$.

We define the linearized residual
\begin{equation}
    \label{eq:residual}
    \delta(z) = \sum_{j=1}^n \frac{y_j\, w_j}{z - z_j} - f(z) \sum_{j=1}^n \frac{w_j}{z - z_j}.
\end{equation}
Given test points $t_1, \ldots, t_m$ in the domain of $f$, we can express the evaluation of $d$ at those points as
\begin{equation}
    \label{eq:residuals}
    \begin{bmatrix} \delta(t_1) \\ \delta(t_2) \\ \vdots \\ \delta(t_m) \end{bmatrix} = \bm{L} \begin{bmatrix} w_1 \\ w_2 \\ \vdots \\ w_n \end{bmatrix},
\end{equation}
where $\bm{L}$ is the $m \times n$ \emph{Loewner matrix} with entries
\begin{equation}
    \label{eq:loewner}
    L_{ij} = - \frac{f(t_i) - y_j}{t_i - z_j}.
\end{equation}

Since the barycentric formula (\ref{eq:barycentric}) is unaffected by a scaling of the weight vector, we can recognize the least-squares minimization of the vector $\bm{\delta}$ as equivalent to finding the least-significant right singular vector of $\bm{L}$. Once the weights are determined, the true residual $f(z) - r(z)$ be evaluated at the test points. The AAA algorithm selects the test point with the largest magnitude of this residual as the next interpolation node.

In the original discrete AAA algorithm~\cite{nakatsukasaAAAAlgorithm2018}, the test points are selected prior to the iteration. In the continuum version~\cite{driscollAAARational2024}, a small initial set of test points is selected to discretize the boundary of the domain, and with each transfer of a test point to the node set, a few additional test points are added along the boundary to either side of the new node. This adaptivity becomes important when $f$ has a singularity very close to or on the boundary of the domain.

\subsection{Greedy Thiele}
\label{sec:thiele}

The Thiele representation of a rational interpolant is~\cite{celisAdaptiveThiele2023,milne-thompsonCalculusFinite1933}
\begin{equation}
    \label{eq:thiele}
    r(z) = d_1 + \cfrac{z - z_1}{d_2 + \cfrac{z - z_2}{d_3 + \cdots}},
\end{equation}
where, given nodes $z_1,\ldots , z_n$ and inductively $d_1,\ldots,d_{n-1}$, the new weight $d_n$ is defined as the result $u_n$ of the iteration
\begin{equation}
    \label{eq:thiele-weight}
    \begin{split}
        u_1 &= f(z_n), \\
        u_{k+1} &= \frac{z_n - z_k}{u_k - d_k} \quad \text{for } k=1,\ldots,n-1.
    \end{split}
\end{equation}
Generically, $r$ interpolates the values $(z_j,f(z_j))$ for $j=1,\ldots, n$ and is a rational function of type $(m,m)$ if $n=2m-1$ or type $(m+1,m)$ if $n=2m$. 

The Thiele representation has had modest computational use due to the inverse difference appearing in~(\ref{eq:thiele-weight}), which is numerically unstable for many node distributions/orderings. However, Salazar~\cite{celisAdaptiveThiele2023,celisNumericalContinued2024} observed that if one follows the structure of AAA by evaluating a candidate $r$ at test points and greedily selecting the worst test point as a new node, the method is often stable. A key attraction of Salazar's method is that it requires only $O(n)$ work to add a new node, comparing favorably to $O(mn^2)$ for AAA with $n$ nodes and $m$ test points.

\subsection{Prescribed poles}
\label{sec:parfrac}

When prior information about the singularity structure of a function is available, it can be useful to use a rational approximation with prescribed poles, yielding a linear problem. Given poles $\zeta_1,\ldots,\zeta_\nu$, consider the partial fractions representation
\begin{equation}
    \label{eq:parfrac}
    r(z) = p(z) + \sum_{j=1}^\nu \frac{c_j}{z - \zeta_j},
\end{equation}
where $q$ is a polynomial. Evaluation of $r$ at the test points $t_1,\ldots,t_m$ yields an overdetermined linear system for the coefficients of $p$ and the residues $c_1,\ldots,c_\nu$.

A potential difficulty with this approach is finding a well-conditioned basis for $p$. On a circle, the monomial basis is ideal, and on an interval, the Chebyshev basis is a natural choice, but on a more general domain, the choice is not as apparent. To resolve this problem, we use the Vandermonde–Arnoldi approach~\cite{BrubeckVandermondeArnoldi2021}, which is a stable method to orthogonalize the monomial basis $\{1,z,z^2,\ldots,z^{N}\}$ as evaluated at the test points $\bm{t}$:
\begin{algorithmic}
    \State $\bm{q}_0 \gets [1,1,\ldots,1]$ 
    \For{$j = 0,\dots,N-1$}
        \State $\bm{v} = \bm{t} \odot \bm{q}_{j}$ \Comment{Hadamard product}
        \For{$k = 0,\dots,j-1$}
            \State $H_{k,j} = \bm{q}_j^* \bm{v} / m$ \Comment{Inner product}
            \State $\bm{v} \gets \bm{v} - H_{k,j} \bm{q}_k$ \Comment{Orthogonalize}
        \EndFor
        \State $H_{j+1,j} = \| \bm{v} \| / \sqrt{m}$ \Comment{Normalize}
        \State $\bm{q}_{j+1} \gets \bm{v} / H_{j+1,j}$ \Comment{Next basis vector}
    \EndFor
\end{algorithmic}
The vectors $\bm{q}_0,\ldots,\bm{q}_N$ form an orthonormal basis for the evaluated monomials, while the upper Hessenberg $\bm{H}$ holds values needed to reconstruct the orthogonalization process at new points. Given coefficients $a_0,\ldots,a_N$ of polynomial $q$ in the orthonormal basis, we evaluate $y=p(z)$ via
\begin{algorithmic}
    \State $y \gets a_0$
    \State $q_0 \gets 1$
    \For{$j = 0,\ldots,N-1$}
        \State $v \gets z q_j$
        \For{$k = 0,\ldots,j-1$}
            \State $v \gets v - H_{k,j} q_k$ 
        \EndFor
        \State $q_{j+1} \gets v / {H}_{j+1,j}$ 
        \State $y \gets y + a_{j+1} q_{j+1}$ 
    \EndFor
\end{algorithmic}

\section{Key data types}
\label{sec:types}

\subsection{DiscretizedPath}

The \textsf{DiscretizedPath} type manages the details of choosing nodes and test points along a \textsf{Curve} or \textsf{Path} from the \textsf{ComplexRegions} package. Both kinds of points are stored in an array that is allocated at the time of construction, with nodes appearing in the first column. The \verb|add_node!| method designates one of the test points as a new node, and additional test points are added to either side of the new node. The \verb|collect| method extends \verb|Base.collect| to return a vector of all nodes, all test points, or all points, depending on the \verb|which| keyword.

\subsection{ArnoldiBasis and ArnoldiPolynomial}

An \textsf{ArnoldiBasis} applies the Arnoldi iteration to a given vector of points, as described in \autoref{sec:parfrac}. The \textsf{ArnoldiPolynomial} type collects an \textsf{ArnoldiBasis} and a vector of coefficients, providing an evaluation method. As a convenience, the backslash operator is overloaded to solve a linear least-squares problem for a function, as in the following:
\begin{lstlisting}[language = Julia, numbers=left, label={lst:exmplg}, caption={Least-squares for an \textsf{ArnoldiBasis}.}]
z = point(Circle(0, 1), range(0, 1, 800))
B = ArnoldiBasis(z, 10)
p = B \ cos    # creates an ArnoldiPolynomial
maximum(abs, p.(z) - cos.(z))
\end{lstlisting}
The result from the above is about $2.78 \times 10^{-7}$. 

\subsection{AbstractRationalFunction and AbstractRationalInterpolant}

The \textsf{AbstractRationalFunction} type is an abstract type providing a common interface for rational functions. The type is parameterized by the floating-point type of the function's values. It provides stubs or universal implementations for the methods shown in \autoref{tab:ratfun_methods}. 

\begin{table}
    \caption{AbstractRationalFunction interface}
    \begin{tabular}{l|l|c}\hline
        Method & Description & Implemented? \\\hline
        \verb|show(io, r)| & Pretty-print details & yes \\
        \verb|r(z)| & Evaluate at a point & no \\
        \verb|eltype(r)| & Floating-point type for values & yes \\
        \verb|isempty(r)| & Is the interpolant empty? & yes \\
        \verb|degrees(r)| & Degree of numerator, denominator & no \\
        \verb|degree(r)| & Degree of the denominator & no \\ 
        \verb|poles(r)| & Poles of the rational function & no \\
        \verb|residues(r)| & Poles and residues & no \\
        \verb|roots(r)| & Roots of the rational function & no \\
        \verb|Res(r, z)| & Residue by trapezoidal rule & yes
    \end{tabular}
    \label{tab:ratfun_methods}
\end{table}

The \textsf{AbstractRationalInterpolant} type is an abstract subtype of \textsf{AbstractRationalFunction} that provides an interface for rational interpolants. It is parameterized by the floating-point type of the interpolation nodes and the floating-point type of the values at the nodes. The type provides stubs or universal implementations for the methods shown in \autoref{tab:ratinterp_methods}.

\begin{table}
    \caption{AbstractRationalInterpolant interface}
    \begin{tabular}{l|l|c}\hline
        Method & Description & Implemented? \\\hline
        \verb|show(io, r)| & Pretty-print details & yes \\
        \verb|nodes(r)| & Interpolation nodes & no \\
        \verb|values(r)| & Values at the nodes & no \\
    \end{tabular}
    \label{tab:ratinterp_methods}
\end{table}

Each following concrete realization of one of these abstract types also implements two signatures of the \texttt{approximate} function that is described in \autoref{sec:approximation}, one for approximation of a function on a continuum domain, and one for fully discrete approximation.

\subsubsection{Barycentric}

\textsf{Barycentric} is a concrete subtype of \textsf{AbstractRationalInterpolant} that implements the barycentric representation of a rational interpolant, as described in \autoref{sec:aaa}. It stores the nodes, values, and weights, and provides implementations for the methods shown in \autoref{tab:ratinterp_methods}. When precision beyond standard \verb|Float64| is detected, the pole computation relies on \texttt{GenericSchur} from \textsf{GenericLinearAlgebra} to perform the necessary generalized eigenvalue computation. T

\subsubsection{Thiele}

\textsf{Thiele} is a concrete subtype of \textsf{AbstractRationalInterpolant} that implements the Thiele form of a rational interpolant, as described in \autoref{sec:thiele}. It stores the nodes, values, and weights, and provides implementations for the methods shown in \autoref{tab:ratinterp_methods}. Currently, however, the implementation of \verb|residues| is based only on trapezoidal integration rather than any special property of the Thiele representation. 

\subsubsection{PartialFractions}

The \textsf{PartialFractions} type is a concrete subtype of \textsf{AbstractRationalFunction} that collects an \textsf{ArnoldiPolynomial} and vectors of poles and residues, representing a rational function with prescribed poles, as given in \autoref{eq:parfrac}. It implements all the methods shown in \autoref{tab:ratfun_methods}. 

\subsection{Approximation}
\label{sec:approximation}

The \textsf{Approximation} type is the primary type that a user interacts with directly. It collects a function, information about its domain, a rational approximant, and an iterative history. The standard way to create an \textsf{Approximation} is via the \verb|approximate| function, which has several different calling signatures. 

\textbf{Interpolation.} For greedy iterations with the barycentric (i.e., AAA) or Thiele representations, \texttt{approximate} takes a \texttt{method} keyword that is either \texttt{Barycentric} (the default) or \texttt{Thiele}. Nodes are added until the error at all test points is less than a tolerance that can be chosen with the keyword \texttt{tol}, or until a stagnation criterion is met for a number of iterations that can be set with the keyword \texttt{stagnation}. The algorithm only returns a result whose poles are all allowed, using a Boolean-valued function that can be set with the \texttt{allowed} keyword. The following signatures are available:
\begin{description}
    \item[\texttt{approximate(f,d)}] for \textsf{Curve} or \textsf{Path} \verb|d| creates a rational interpolant of the function \verb|f| on the domain \verb|d| using a continuum approach. No poles in \verb|d| are allowed.
    \item[\texttt{approximate(f,d)}] for \textsf{SimplyConnectedRegion} \verb|d| creates a rational interpolant of the function \verb|f| on the domain \verb|d| using a continuum approach. The nodes are chosen on the boundary of \verb|d|, and no poles are allowed in \verb|d|, which may be an interior or exterior region.
    \item[\texttt{approximate(f,z)}] for a vector of numbers \verb|z| creates a rational interpolant of the function \verb|f| in a discrete approach, using \verb|z| as the fixed set of test points and potential nodes. By default, all poles are allowed.
    \item[\texttt{approximate(y,z)}] for vectors \verb|y| and \verb|z| uses the entries of \verb|y| as the values to be taken at the entries of \verb|z|. This method returns a subtype of \textsf{AbstractRationalInterpolant}, as there is presumed to be no information about an originating function.
\end{description}

Because both the AAA and Salazar iterations rely on pairwise differences between nodes, test points, and function values in various combinations, care was taken to reuse such information between iterations. The representation-specific references in the iteration are the representation constructor and the methods \verb|add_nodes!| and \verb|update_test_values!|, which the representations must implement. A first call to \verb|update_test_values!| allocates working space. Subsequent calls update the working arrays at any new test points, updates the weights of the rational interpolant (in the barycentric case only), and evaluates the interpolant at all test points. Once a new node is chosen, a call to \verb|add_nodes!| updates the node and value properties of the interpolant (which, in the Thiele case, also updates the weights), allocates space for new weights, and updates differencing information between the new node(s) and the existing test points. 

\textbf{Least squares.} For a least-squares approximation with prescribed poles using the partial fraction representation, a \texttt{degree} keyword to \texttt{approximate} can be used to specify the degree of the polynomial part. The following signatures are available:
\begin{description}
    \item[\texttt{approximate(f,d,zeta)}] for \textsf{Curve} or \textsf{Path} \verb|d| creates an approximation of the function \texttt{f} on the domain \verb|d| with the prescribed poles in the vector \verb|zeta|. To discretize \verb|d|, equally spaced nodes are placed (whose number is selectable by the \verb|init| keyword), and test points are promoted to nodes iteratively until the distance between adjacent nodes is no greater than one-half the distance to the nearest pole. 
    \item[\texttt{approximate(f,z,zeta)}] for a vector of numbers \verb|z| uses the entries of \verb|z| as the discretization of the domain.
    \item[\texttt{approximate(y,z,zeta)}] for vectors \verb|y| and \verb|z| uses the entries of \verb|y| as the values to take at the entries of \verb|z|. This method returns a subtype of \textsf{AbstractRationalFunction}, as there is presumed to be no information about an originating function.
\end{description}

\section{Examples}

\begin{figure}
\centering
\includegraphics[width=0.75\columnwidth]{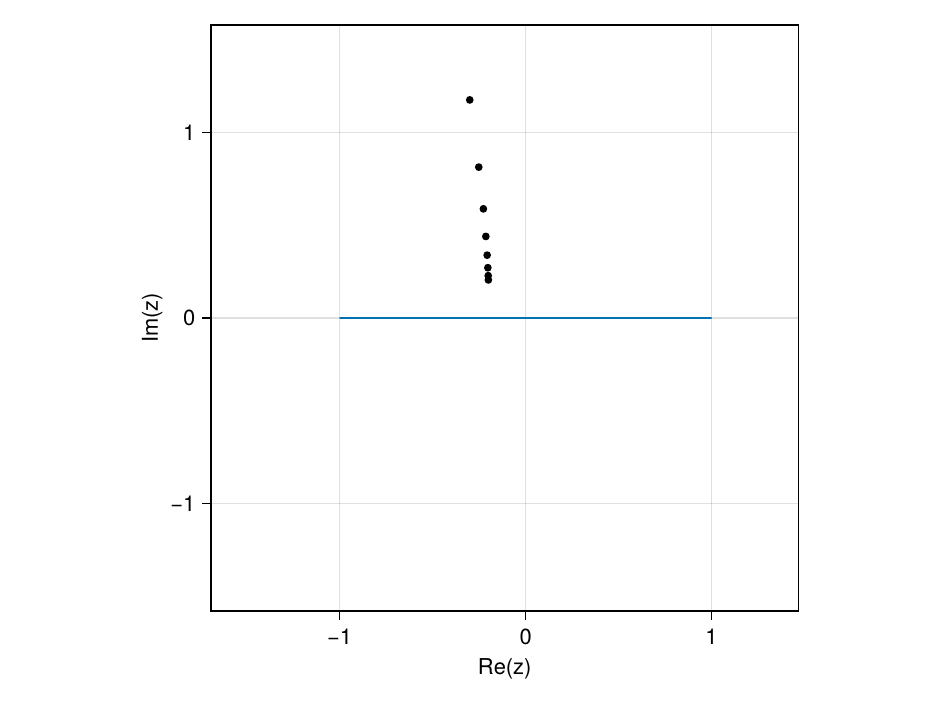}
\caption{Pole plot of the AAA rational interpolant of $f(z) = \log(1 + i + 5iz)$ on the unit interval. The poles cluster near the branch point at $(i-1)/5$.}
\label{fig:poleplot1}
\end{figure}

The package defines \verb|unit_interval| as a \textsf{Segment} from $-1$ to $1$. We can approximate the function $f(z) = \log(1 + i + 5iz)$ on this interval as follows:
\begin{lstlisting}[language = Julia, caption={Continuum AAA on the unit interval.}, label={lst:log}]
using RationalFunctionApproximation, CairoMakie
f = z -> log(1 + 1im + 5im*z)
r = approximate(f, unit_interval)
poleplot(r)
\end{lstlisting}
The result is a \textsf{Barycentric} interpolant of type $(12,12)$ and the plot in \autoref{fig:poleplot1}. The poles of this interpolant cluster near the branch point of $f$ at $(i-1)/5$. We can check the size of the error via \verb|check(r)|, which returns vectors of test points and error values and prints a message with the maximum error; here, that maximum is about $1.6\times 10^{-13}$. 

\begin{figure}
\centering
\includegraphics[width=0.75\columnwidth]{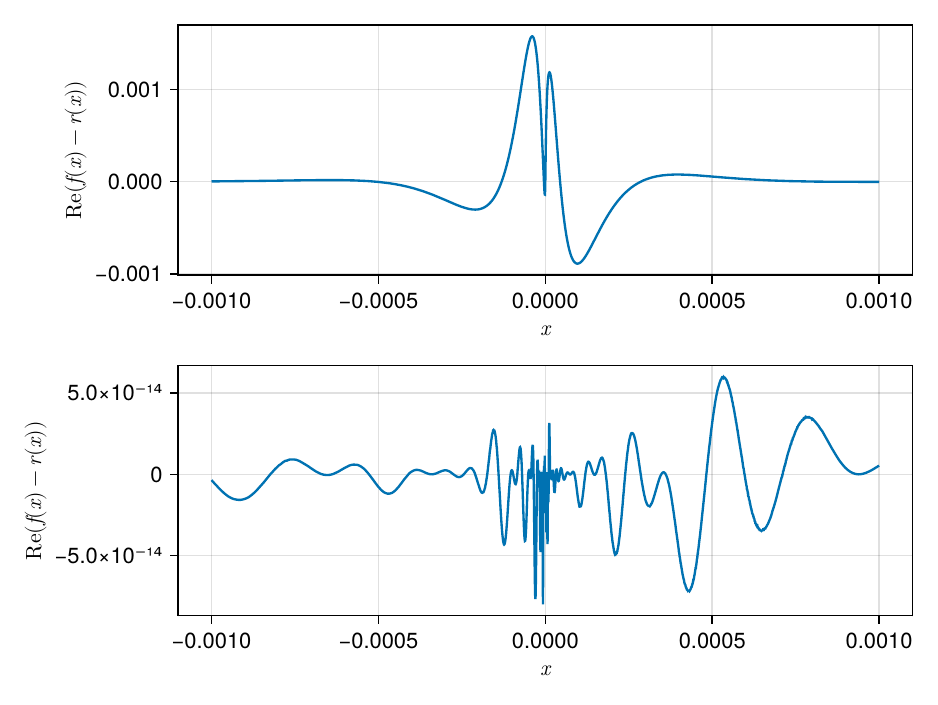}
\caption{Top: In the discrete AAA iteration, accuracy close to a singularity will suffer if the nodes have insufficient resolution. Bottom: Continuum AAA adaptively chooses nodes and test points in order to resolve the singularity.}
\label{fig:continuum}
\end{figure}

The usefulness of the continuum methods is apparent when singularities of $f$ creep close to the domain. If, for instance, we approximate $\sqrt{x + 10^{-6}i}$ on equally spaced points in $[1,1]$, we will miss important details at the origin:
\begin{lstlisting}[language = Julia, caption={Discrete AAA for $\sqrt{x + 10^{-6}i}$.}]
f = z -> sqrt(z + 1e-6im)
x = range(-1, 1, 1001)
r = approximate(f, x)
lines(-0.001..0.001, x->real(f(x) - r(x)))
\end{lstlisting}
Although the discrete iteration finishes successfully, with a reported maximum error of $4.5\times 10^{-14}$ at the entries of \verb|x|, we see in \autoref{fig:continuum} that the error near the origin exceeds $10^{-3}$ at other locations. With advance knowledge of the singularity structure, a discretization could be chosen manually that would properly resolve the singularity. But the continuum approach offers this resolution automatically:
\begin{lstlisting}[language = Julia, caption={Continuum AAA for $\sqrt{x + 10^{-6}i}$.}]
r = approximate(f, unit_interval)
lines(-0.001..0.001, x->real(f(x) - r(x)))
\end{lstlisting}
Now, as confirmed in \autoref{fig:continuum}, the error is less than $10^{-13}$ everywhere on the interval. 

\begin{figure}
\centering
\includegraphics[width=0.75\columnwidth]{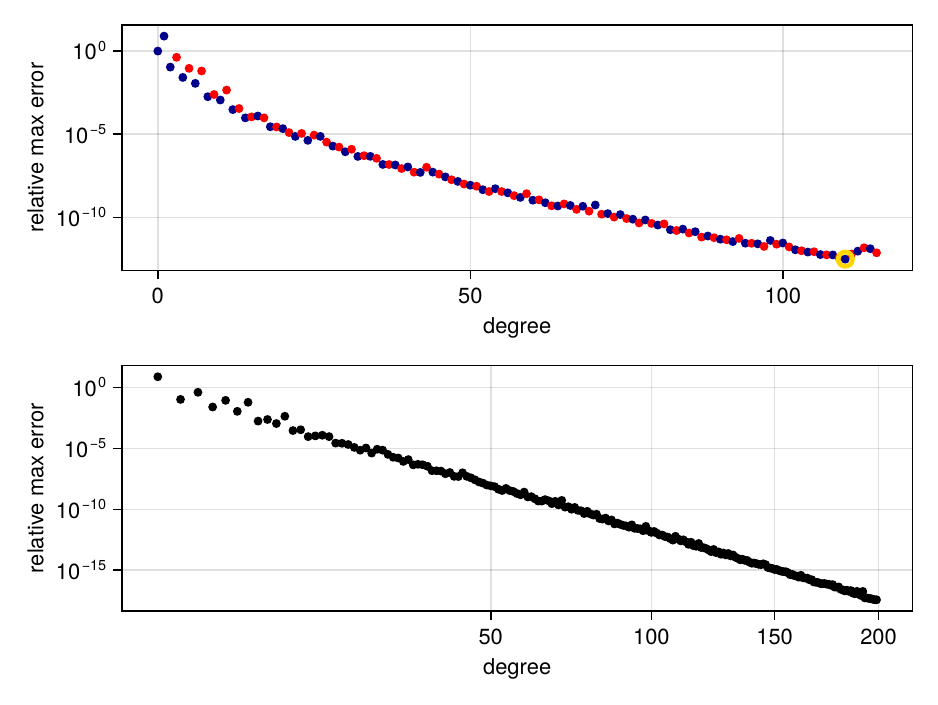}
\caption{Convergence of continuum AAA to the function $|x|$ over $[-1,1]$. Top: Red dots indicate an interpolant having a pole in the interval, while blue dots indicate no disallowed poles. The gold halo shows the interpolant returned as the result, chosen when the iteration stagnated. Bottom: Using extended precision, the iteration continues successfully (without pole checking); the square-root scale on the degree axis illustrates root-exponential convergence.}
\label{fig:absconverge}
\end{figure}

\begin{figure}
\centering
\includegraphics[width=0.75\columnwidth]{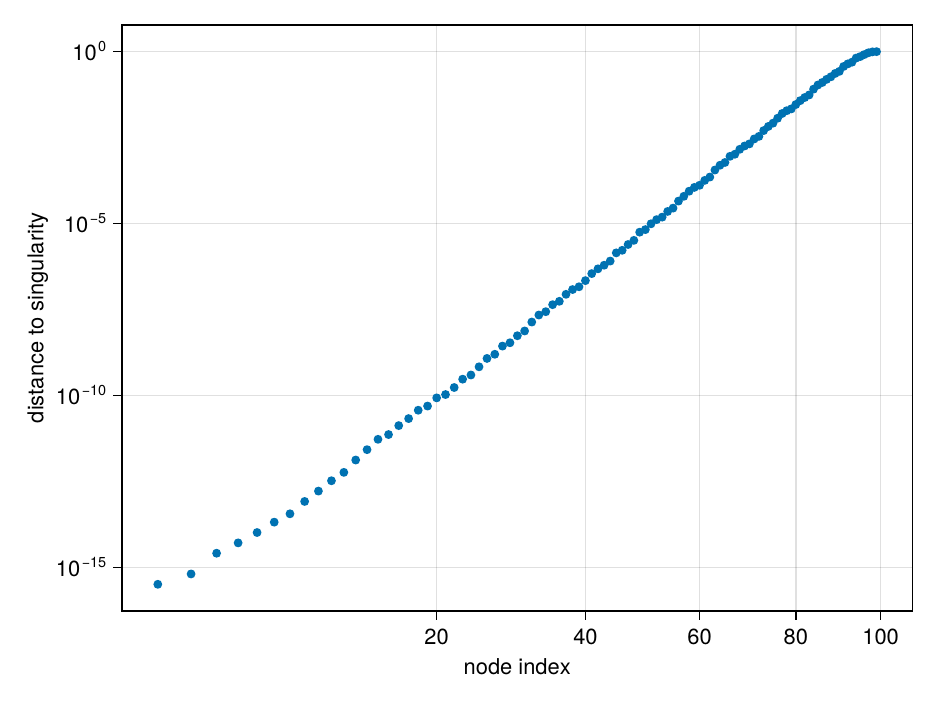}
\caption{Positive nodes selected by continuum AAA for the function $|x|$ over $[-1,1]$. They cluster near the singularity at $x=0$.}
\label{fig:absnodes}
\end{figure}

The convergence of AAA for Newman's example of $f(x)=|x|$ is interesting:
\begin{lstlisting}[language = Julia, caption={Continuum AAA for $|x|$.}, label={lst:abs}]
r = approximate(abs, unit_interval)
convergenceplot(r)
\end{lstlisting}
As seen in the upper plot of \autoref{fig:absconverge}, the iteration was halted when convergence stalled a bit short of machine roundoff. In order to continue past the plateau, we can use extended precision supplied by the \textsf{DoubleFloats} package. Since the eigenvalue computation needed to check for disallowed poles is much slower in extended precision, we disable the checks with the keyword \verb|allowed=true|. 
\begin{lstlisting}[language = Julia, caption={Continuum AAA in extended precision.}, label={lst:absquad}]
using DoubleFloats, ComplexRegions
uiquad = Segment{Double64}(-1, 1)
r = approximate(abs, uiquad; 
        allowed=true, max_iter=200)
\end{lstlisting}
This computation took 30.5 seconds to compute on an Apple M4 Pro processor, compared to 0.23 seconds in \autoref{lst:abs}. The lower plot of \autoref{fig:absconverge} shows that the convergence continues smoothly at a root-exponential rate. \autoref{fig:absnodes} illustrates that the positive node locations cluster with a tapered exponential pattern at the singularity, much as the poles do~\cite{trefethenExponentialNode2021}. 

\begin{figure}
\centering
\includegraphics[width=0.75\columnwidth]{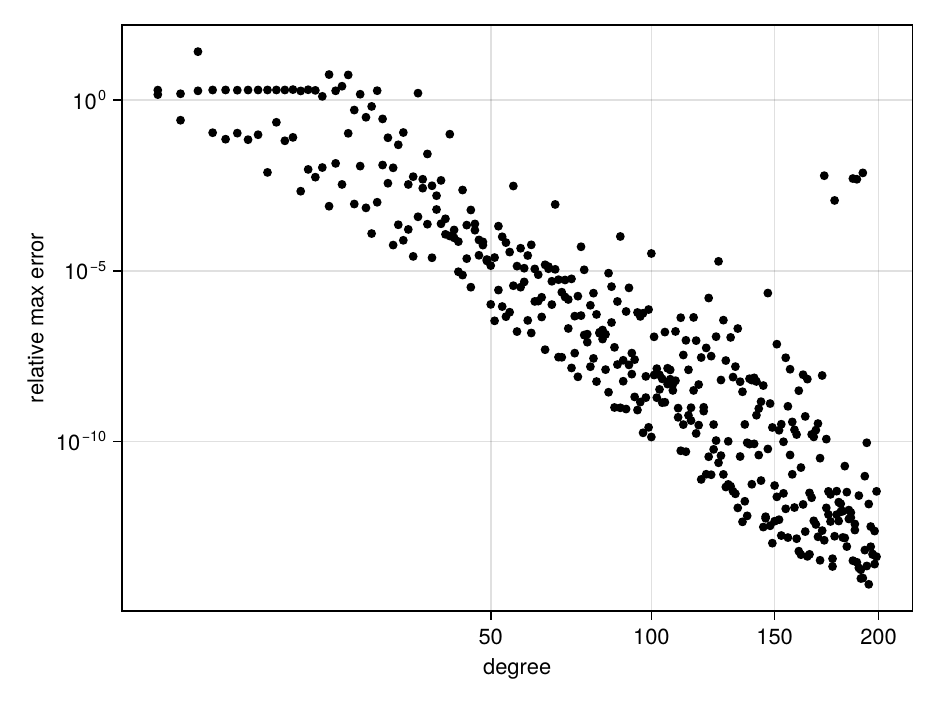}
\caption{Convergence of continuum Thiele to the function $|x|$ over $[-1,1]$. The smallest errors are similar to those in \autoref{fig:absconverge} for AAA, but the convergence is far less smooth.}
\label{fig:absthiele}
\end{figure}

Any of the above examples can be computed using the greedy Thiele iteration in place of AAA by adding \verb|method=Thiele| to the \verb|approximate| call. For example, we can use it in extended precision for $|x|$:In order to exploit its speed for \autoref{lst:abs}, we use \verb|allowed=true| to disable pole checking and add \verb|stagnation=20| to make it persist through plateaus. We also double the number of iterations, since the degree is roughly half the number of nodes:
\begin{lstlisting}[language = Julia, caption={Continuum Thiele for $|x|$.}]
approximate(abs, uiquad; allowed=true, 
        method=Thiele, max_iter=400, stagnation=20)
\end{lstlisting}
The above code takes just 1.69 seconds to reach a similar degree and accuracy as the AAA result. This speedup is dramatic because AAA relies on generic implementations of the SVD for \textsf{DoubleFloats}, while Thiele requires no linear algebra at all. However, as seen in \autoref{fig:absthiele}, the convergence of the Thiele iteration is far less smooth here, which is why \verb|stagnation=20| was used to make it continue past apparent plateaus.

\begin{figure}
\centering
\includegraphics[width=0.75\columnwidth]{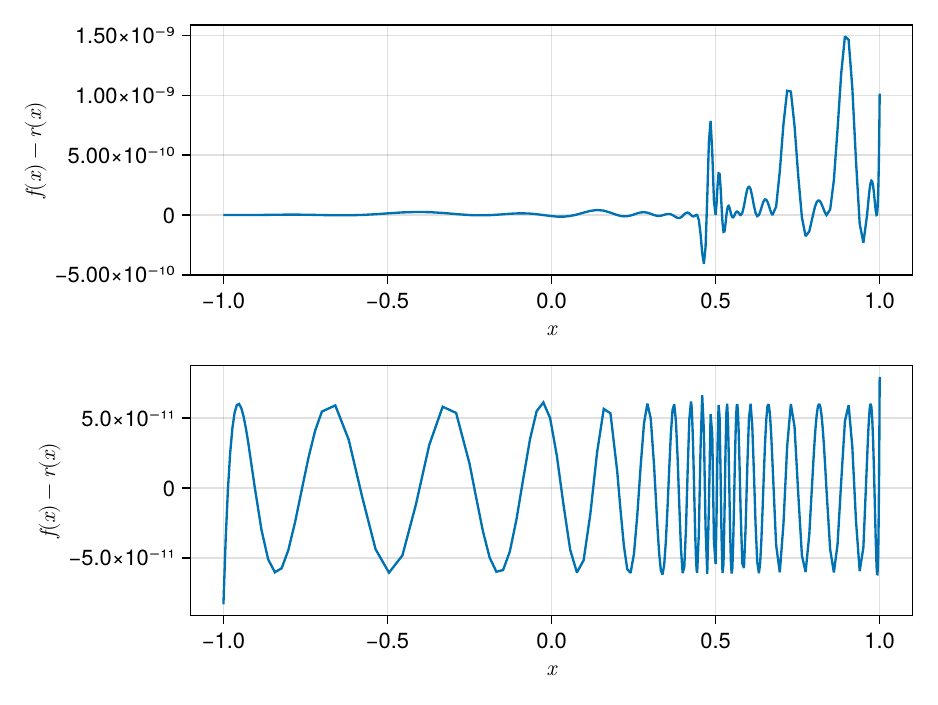}
\caption{Top: Nonuniform error is produced by the AAA iteration, which uses a least-squares criterion. Bottom: After Lawson-style iterations, the error can be made nearly uniform across the interval.}
\label{fig:minimax}
\end{figure}

An advantage of the AAA approach is that it can be modified to construct near-best approximations in the max-norm sense~\cite{nakatsukasaAlgorithmReal2020}, which classically is of greater interest than the 2-norm. For instance, the error for $f(x)=|x-0.5 + 0.05i|$ is rather nonuniform:
\begin{lstlisting}[language = Julia, caption={Continuum AAA for $|x-0.5 + 0.05i|$.}]
f = z -> abs(z - 0.5 + 0.05im)
r = approximate(f, unit_interval, max_iter=20)
lines(check(r)...)
\end{lstlisting}
As seen on the top of \autoref{fig:minimax}, the error is larger on the right end of the interval than on the left by orders of magnitude. The \verb|minimax| function uses iteratively weighted norms to place more emphasis where the error is largest, resulting in near-uniform error, as shown in the bottom of the figure.
\begin{lstlisting}[language = Julia, caption={Minimax AAA for $|x-0.5 + 0.05i|$.}]
r_inf = minimax(r, 20);
lines(check(r_inf)...)
\end{lstlisting}

\begin{figure}
\centering
\includegraphics[width=0.75\columnwidth]{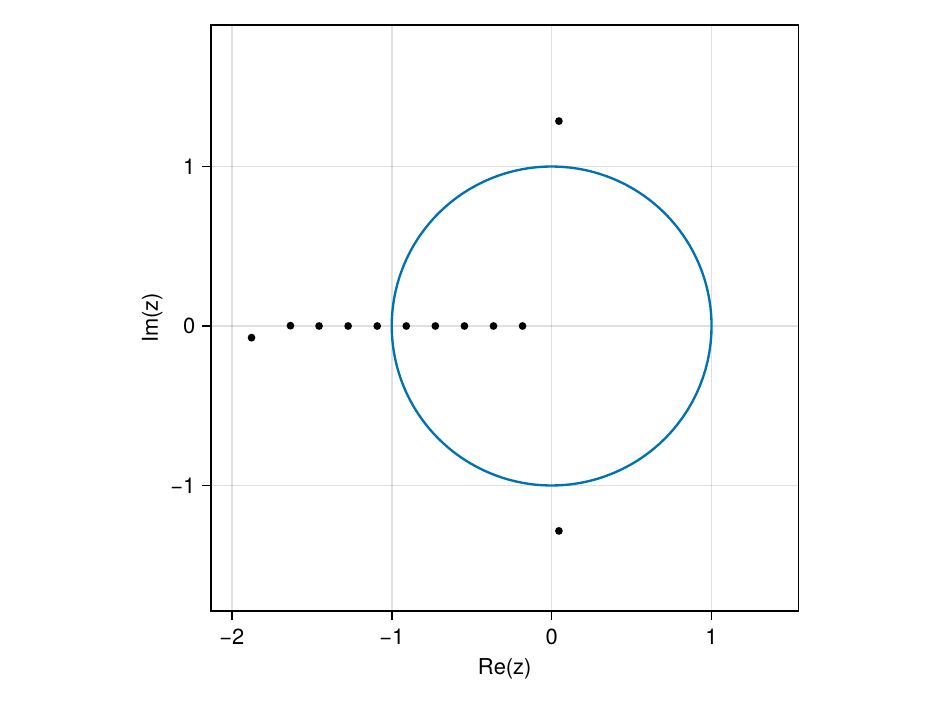}
\caption{Poles of a rational interpolant of $f(z) = 1/\zeta(11z)$ on the unit circle.}
\label{fig:zetapoles}
\end{figure}

\begin{figure}
\centering
\includegraphics[width=0.75\columnwidth]{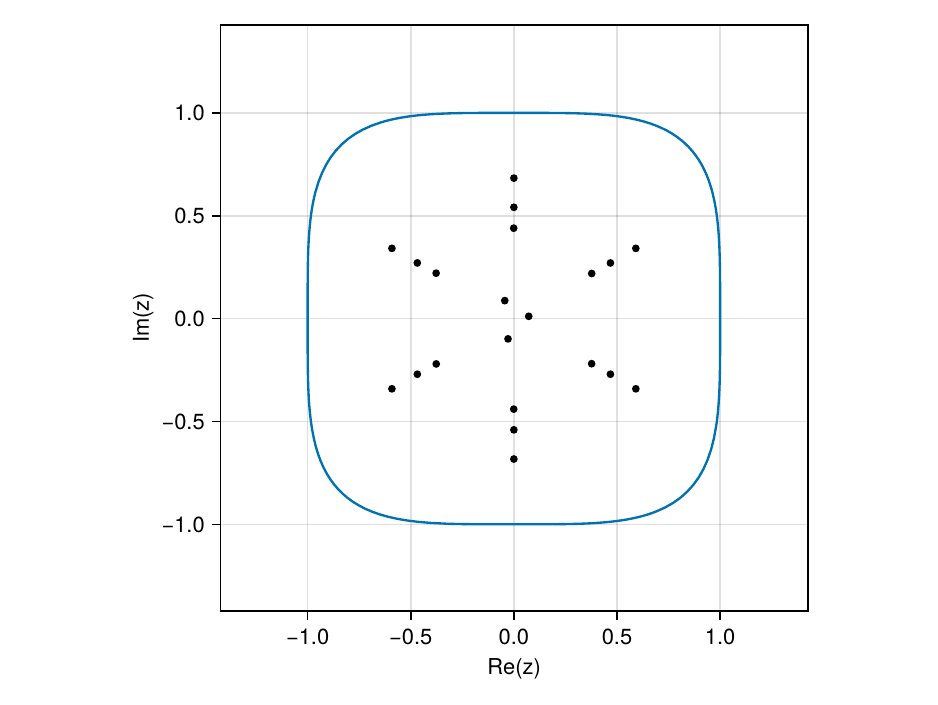}
\caption{Poles of a rational approximation of $\coth(z^{-3})$ in the exterior of a squircle.}
\label{fig:squircle}
\end{figure}

Other domains of analyticity may easily be specified. One predefined domain is the unit circle: 
\begin{lstlisting}[language = Julia, caption={Approximation on the unit circle.}]
using SpecialFunctions: zeta
f = z -> 1/ zeta(11z)
r = approximate(f, unit_circle)
poleplot(r)
\end{lstlisting}
The poles of the resulting interpolant of degree 21 are shown in \autoref{fig:zetapoles}. The maximum error in the approximation on the circle is about $1.2\times 10^{-12}$. The \textsf{ComplexRegions} package defines a \textsf{Shapes} module with other predefined curves, as well as functions \verb|interior| and \verb|exterior| that can be used to create \textsf{SimplyConnectedRegion} objects.
\begin{lstlisting}[language = Julia, caption={Approximation on an exterior domain.}]
f = z -> coth(1/z^3)
r = approximate(f, exterior(Shapes.squircle))
poleplot(r)
\end{lstlisting}
The resulting plot is given in \autoref{fig:squircle}.

If prior knowledge of the singularity structure is available, we can use the partial fraction representation and solve a linear problem without iteration. For example, suppose we return to the log function in \autoref{lst:log} and extract the poles from a rational interpolant. We can use the same poles to approximate a different function with the same branch point:
\begin{lstlisting}[language = Julia, caption={Partial fraction approximation of $\log(1 + i + 5iz)$.}]
f = z -> log(1 + 1im + 5im*z)
r = approximate(f, unit_interval)
ζ = poles(r)
g = z -> sqrt(1 + 1im + 5im*z)
s = approximate(g, unit_interval, ζ; degree=20)
\end{lstlisting}
The resulting approximation reports an error of size $5.3\times{10}^{-10}$ at its test points. 

\section{Future work}

It should be possible to implement operations on the rational function types to enable basic arithmetic, composition, differentiation, and integration, in the style of Chebfun~\cite{driscollChebfunGuide2014a} and \textsf{ApproxFun}~\cite{olverPracticalFramework2014}. A systematic comparison of the AAA and Thiele methods is also needed, as well as a more complete understanding of the stability of the Thiele method. Least-squares approximations with prescribed poles have been fruitful in the solution of Laplace and related PDEs~\cite{costaAAAleastSquares2023,gopalSolvingLaplace2019,gopalRepresentationConformal2019}, which requires imposing conditions other than direct approximation of a known function.

\section{Acknowledgments}
I thank Nick Trefethen and Oliver Salazar Celis for their valuable feedback on a draft of this paper.


\bibliographystyle{juliacon}
\bibliography{ref.bib}

\end{document}